\documentclass{article}

\usepackage{microtype}
\usepackage{verbatim}
\usepackage{graphicx}
\usepackage{rotating}

\usepackage{caption}
\usepackage{subcaption}

\usepackage{url}
\usepackage{hyperref}

\usepackage{pslatex}

\usepackage{datetime}
\newdateformat{versiondate}{%
\THEMONTH\THEDAY}

\usepackage{pgf,tikz,environ}
\usetikzlibrary{arrows}
\usetikzlibrary{cd}

\usepackage{xpatch} 

\usepackage{amsthm}
\usepackage[fleqn]{amsmath} 
\usepackage{amsfonts}
\usepackage{amssymb}

\newtheoremstyle{zoltanstyle}
  {1em} 
  {\topsep} 
  {} 
  {} 
  {\bfseries} 
  {.} 
  {.5em} 
  {} 

\theoremstyle{zoltanstyle}
\swapnumbers
\xpatchcmd\swappedhead{~}{.~}{}{}
\newtheorem{body}{}
\numberwithin{body}{section}

\newtheorem{corollary}[body]{Corollary}
\newtheorem{definition}[body]{Definition}

\newtheorem{lemma}[body]{Lemma}

\newtheorem{proposition}[body]{Proposition}

\newtheorem{theorem}[body]{Theorem}

\expandafter\let\expandafter\oldproof\csname\string\proof\endcsname
\let\oldendproof\endproof

\renewenvironment{proof}[1][\proofname]{%
  \oldproof[\normalfont \bfseries #1.]%
}{\oldendproof}

\usepackage{enumitem}
\usepackage{mathptmx}
\usepackage{xspace}

\usepackage[T1]{fontenc}
\usepackage[utf8]{inputenc}

\newcommand{\SetComp}[2]{\left\{ {#1}\:\middle|\:{#2} \right\}}   
\newcommand{\Bool}{\mathbf{2}} 
\newcommand{\jequiv}{=} 
\newcommand{\pet}{\mathrm{pe}_\mathbf{1}} 
\newcommand{\peb}{\mathrm{pe}_\mathbf{0}} 

\title{Apartness relations between propositions}
\author{Zoltan A. Kocsis \thanks{Computer Science and Engineering, University of New South Wales, Kensington NSW 2052, Australia.}}
\date{\today}

\usepackage[normalem]{ulem} 

\begin{document}

\maketitle

\begin{abstract}
We classify all apartness relations definable in propositional logics extending intuitionistic logic using Heyting algebra semantics. We show that every Heyting algebra which contains a non-trivial apartness term satisfies the weak law of excluded middle, and every Heyting algebra which contains a tight apartness term is in fact a Boolean algebra. This answers a question of E.~Rijke regarding the correct notion of apartness for propositions, and yields a short classification of apartness terms that can occur in a Heyting algebra. We also show that Martin-Löf Type Theory is not able to construct non-trivial apartness relations between propositions.
\end{abstract}

\section{Introduction}\label{sec1}

\begin{definition}\label{def:apartness-rel}
Given a set $S$, we call an irreflexive, symmetric binary relation $\# \subseteq S \times S$ an \textit{apartness} relation if the following condition (cotransitivity) holds for all $x,y,z \in S$: if $x \# y$ then either $x \# z$ or $z \# y$. We call the apartness relation \textit{tight} if the following condition is satisfied as well: $\text{if } x \# y \text{ does not hold, then } x = y.$
\end{definition}

\begin{body}
The complement of an apartness relation clearly yields an equivalence relation on the set $S$. When working in classical logic, every apartness relation itself arises as the complement of the equivalence relation obtained this way: this completely reduces the study of apartness relations to the familiar subject of equivalence relations. Moreover, classical foundational theories trivialize the notion of tight apartness altogether, as the only tight apartness arises as the complement of the equality relation itself.
\end{body}

\begin{body}
When we work in an intuitionistic foundational theory, the complement of the equality relation on a set $S$ may fail to constitute an apartness relation, and a carefully chosen apartness often has more constructive content than (in)equality itself. In constructive analysis, apartness relations are encountered  frequently; for example, in a metric space, defining $x$ to be apart from $y$ if $d(x,y)>0$ provides an apartness relation. Likewise, one of the two main schools of constructive algebra uses apartness relations extensively, studying structures equipped with a tight apartness and operations that respect them, while taking quotients via anti-subalgebras and anti-ideals, which behave better constructively than the classical notion of (complements) of subalgebras and ideals~\cite{tvdConst}.
\end{body}

\begin{body}
In this article we investigate apartness relations between intuitionistic propositions. Motivation for the study of such relations comes from multiple sources. For example, in Homotopy Type Theory, the propositions-as-some-types paradigm, which identifies propositions with subsingleton types, supplies many interesting irreflexive and symmetric relations on the lattice of propositions; among them questions about the relation $P \not\sim Q$ asserting the contractibility of the coproduct type $P + Q$ for propositions $P,Q$. One would like to know the circumstances in which this and other similar relations satisfy tightness/cotransitivity. One can show that (as long as $P,Q$ have at most one inhabitant) $P \not\sim Q$ holds precisely if the type $(P \times \neg Q) + (\neg P \times Q)$ has an inhabitant. As we shall see below, this means that $\not\sim$ satisfies cotransitivity precisely if the Law of Excluded Middle holds for all propositions.
\end{body}

\begin{body}
Here, by means of Heyting algebra semantics, we classify all apartness relations definable in propositional logics extending intuitionistic logic. We show that every Heyting algebra which contains a non-trivial apartness term satisfies the weak law of excluded middle (Theorem~\ref{thm:characterization}), and every Heyting algebra which contains a tight apartness term is in fact a Boolean algebra (Theorem~\ref{thm:boolean}). This answers a question of E.~Rijke~\cite{erijke} on the correct notion of apartness for propositions, and yields a short classification of apartness terms that can occur in a Heyting algebra (Theorem~\ref{thm:classification}). We further observe that Martin-Löf Type Theory does not define any non-trivial apartness relations between propositions (Corollary~\ref{cor:mlttnondef}) under the propositions-as-types or propositions-as-some-types paradigms. Finally, with some additional work to eliminate appeals to semantic observations (Lemma~\ref{lemma:typed-reducts} and Lemma~\ref{lemma:typed-reducts-converse}), we show that the classification results obtained for intuitionistic logic hold in Martin-Löf Type Theory as long as one uses the latter paradigm.
\end{body}

\begin{body}
The relation $\#_2$, used in Section~\ref{sec:typetheory}, was first introduced by M.~Escardó~\cite{escardoTypeTop} as a way of characterizing \textit{totally separated} types: a type is totally separated precisely if this apartness relation is tight. The most recent revisions of that work contain a formalized, type-theoretic analogue of Proposition~\ref{prop:lem-rijke1} as well, under the name \texttt{$\#\Omega$-cotran-taboo}.
\end{body}

\begin{body}
T.~de~Jong~\cite{deJongApartness} investigated apartness relations on continuous directed-complete partial orders (DCPOs), and showed that in a constructive setting, nontrivial DCPOs cannot admit nontrivial apartness relations. More precisely, if $x \leq y$ and $x \# y$ both hold in such a DCPO, then the weak law of excluded middle holds for any proposition in the metatheory. This bears some resemblance to Theorem~\ref{thm:characterization}, which also concludes the weak law of excluded middle from an apartness assumption, but in the algebra itself instead of the metatheory.
\end{body}

\begin{body}
According to A.~S.~Troelstra and D.~van~Dalen~\cite{tvdConst} (8.1 in Chapter 5), Theorem~\ref{thm:stability} on the double-negation stability of equality can be traced back to L.~E.~J.~Brouwer~\cite{brouwerStab}; a short, general proof is given in the second volume~\cite{tvdConst2} of the same work (1.2 in Chapter 8).
\end{body}

\subsection{Notation and Preliminaries}

\begin{body}
Recall that the standard semantics of formulas of the classical propositional calculus interprets each formula $\varphi$ of classical propositional logic as a corresponding \textit{term} $\varphi'$ in the language of Boolean algebras. Soundness and completeness of the semantics shows that one can prove a given formula $\varphi$ in classical propositional logic precisely if the equation $\varphi'= \top$ holds for the corresponding term $\varphi'$ in \textit{all} Boolean algebras.
\end{body}  

\begin{body}
A \textit{Heyting algebra} $\langle H, \wedge, \vee, \rightarrow, \bot, \top \rangle$ is a bounded distributive lattice where the set $\SetComp{c \in H}{a \wedge c \leq b}$ has a distinguished greatest element (denoted $a \rightarrow b$) for any pair of elements $a, b \in H$.
\end{body}

\begin{body}
In what follows, let $\mathcal{L}$ denote the first-order language of Heyting algebras. As usual, we use $\neg x$ as shorthand for $x \rightarrow \bot$ in this language, and we abbreviate $(x \rightarrow y) \wedge (y \rightarrow x)$ as $x \leftrightarrow y$. Moreover, we write $a \leq b$ if all of the three (equivalent) conditions $a \wedge b = a$, $a \vee b = b$, $a \rightarrow b = \top$ hold.
\end{body}

\begin{body}
Similarly to the standard semantics of classical propositional logic, we can identify the propositional formulae of intuitionistic logic IPC with first-order \textit{terms} in the language $\mathcal{L}$. We get the soundness-completeness result that one can prove an intuitionistic propositional formula $\varphi$ in IPC (which we denote as $\vdash \varphi$) precisely if the equation $\varphi' = \top$ holds for the corresponding $\mathcal{L}$-term $\varphi'$ in every Heyting algebra $H$ (which we denote as $H \models \varphi' = \top$). For an overview of Heyting algebra semantics for intuitionistic propositional logic specifically, consult \textit{A semantic hierarchy for Intuitionistic Logic}~\cite{gbSemantic}.
\end{body}

\begin{body}
Section~\ref{sec:typetheory} concerns Martin-Löf Type Theory. That section begins with a brief outline of the relevant syntax of Type Theory; where not otherwise noted, we stick to the conventions set out in version \texttt{first-edition-1351-g99f4de9} of \textit{The HoTT Book} \cite{hottbook}.
\end{body}

\subsection{Heyting algebras}

\begin{body}
We call the equation $x \vee \neg x = \top$ the \textit{law of excluded middle} and call the equation $\neg x \vee \neg\neg x = \top$ the \textit{weak law of excluded middle}. A Heyting algebra $H$ is a Boolean algebra precisely if the law of excluded middle (or equivalently the equation $\neg \neg x = x$) holds for all of its elements, and in that case $x \rightarrow y = \neg x \vee y$ also holds for all $x,y \in H$.
\end{body}

{ \allowdisplaybreaks
\begin{definition}
Consider an IPC formula $\varphi$ and a proposition letter $P$. As usual, we define the \textit{substitution operator} $[P/\varphi]$ on IPC formulas recursively, using the following equations:
\begin{align*}
    & [P/\varphi](P) = \varphi, & \text{} \\
    & [P/\varphi](X) = X, & \text{for atomic $X \neq P$} \\
    & [P/\varphi](\psi_1 \wedge \psi_2) = [P/\varphi](\psi_1) \wedge [P/\varphi](\psi_2), & \text{} \\
    & [P/\varphi](\psi_1 \rightarrow \psi_2) = [P/\varphi](\psi_1) \rightarrow [P/\varphi](\psi_2), & \text{} \\
    & [P/\varphi](\psi_1 \vee \psi_2) = [P/\varphi](\psi_1) \vee [P/\varphi](\psi_2) & \text{for any $\psi_1$ and $\psi_2$.}
\end{align*}
As customary, we write the substitution operator in post-fix notation, as $\psi[P/\varphi]$, instead of $[P/\varphi](\psi)$ in the remainder of this article.
\end{definition}
}

\begin{definition}
Consider a set of IPC formulas $L$. We call $L$ (the set of theorems of) a \textit{superintuitionistic logic} if the following hold for all intuitionistic formulas $\varphi, \psi$ and proposition letters $P$:
\begin{enumerate}
    \item If $\vdash \varphi$, then $\varphi \in L$;
    \item if $\varphi \in L$ and $\varphi \rightarrow \psi \in L$, then $\psi \in L$;
    \item if $\varphi \in L$, then $\varphi[P/\psi] \in L$.
\end{enumerate}
\end{definition}

\begin{definition}
Take a superintuitionistic logic $L$ and an IPC formula $P \# Q$ in two proposition letters $P,Q$. Given IPC formulas $\varphi, \psi, \theta$, let $\varphi \# \psi$ denote the result of the substitution $(P \# Q)[P/\varphi][Q/\psi]$. We call $P \# Q$ an \textit{apartness relation between propositions} in the logic $L$ if it satisfies the following:
\begin{itemize}
    \item \textbf{Irreflexive:} $(\varphi \# \varphi) \in L$,
    \item \textbf{Symmetric:} $(\varphi \# \psi) \rightarrow (\psi \# \varphi) \in L$,
    \item \textbf{Cotransitive}: $(\varphi \# \psi) \rightarrow ((\varphi \# \theta) \vee (\theta \# \psi)) \in L$.
\end{itemize}
We call the apartness relation between propositions \textit{tight} if furthermore the formula $\neg(\varphi \# \psi) \rightarrow (\varphi \leftrightarrow \psi)$ belongs to $L$ as well.
\end{definition}

\section{Apartness terms}

\begin{body}
We wish to classify all apartness relations between propositions in superintuitionistic logics. We will accomplish this by appealing to results on free Heyting algebras, most importantly Theorem~\ref{thm:rieger-nishimura}. Definition~\ref{def:apartness-term} provides an analogue of apartness relations for Heyting algebras.
\end{body}

\begin{definition}\label{def:apartness-term}
Consider an $\mathcal{L}$-term $x \# y$ with two free variables $x,y$. We say that $x \# y$ constitutes an \textit{apartness term} in the Heyting algebra $H$ if the following three hold for all $x,y,z \in H$:
\begin{itemize}
    \item \textbf{Irreflexive}: $x \# x = \bot$;
    \item \textbf{Symmetric}: $x \# y = y \# x$;
    \item \textbf{Cotransitive}: $x \# y \leq (x \# z) \vee (z \# y)$.
\end{itemize}
In addition, we call the term $x\# y$ \textit{trivial} if $x \# y = \bot$, and \textit{tight} if $\neg (x \# y) \wedge x \leq y$.
\end{definition}

\begin{body}
Keep in mind that the notion of apartness term does not coincide with the notion of an \textit{apartness relation on a Heyting algebra}: one can equip any Heyting algebra straightforwardly with a (non-classically not necessarily tight) apartness relation just by taking the negation of equality, but we will shortly prove that not all Heyting algebras admit an apartness \textit{term}.
\end{body}

\begin{proposition}\label{prop:lindenbaum}
A superintuitionistic logic $L$ has an apartness (resp. tight apartness) relation between its propositions precisely if its Lindenbaum algebra $L'$ (which is a Heyting algebra) has an apartness (resp. tight apartness) term.
\end{proposition}

\begin{body}
Rijke~\cite{erijke} proposed the following formulas as ``natural candidates'' for non-trivial apartness relations between propositions:
\begin{enumerate}
    \item the formula $(P \wedge \neg Q) \vee (\neg P \wedge Q)$;
    \item the formula $\neg (P \leftrightarrow Q)$.
\end{enumerate}
\end{body}

\begin{body}
We begin by showing that none of Rijke's candidates are apartness relations between propositions in plain IPC. We then give necessary and sufficient conditions on a superintuitionistic logic $L$ which ensure that Rijke's candidates give rise to apartness relations between propositions. Eventually, we will show that, up to logical equivalence, the second of these is in fact the \textit{only} non-trivial candidate in every superintuitionistic logic. In accordance with Proposition~\ref{prop:lindenbaum}, from here on out we pass between apartness terms in Heyting algebras and apartness relations in superintuitionistic logics without further commentary.
\end{body}

\begin{proposition}\label{prop:lem-rijke1}
Rijke's first candidate, the term $(x \wedge \neg y) \vee (\neg x \wedge y)$, constitutes an apartness term on a Heyting algebra $H$ precisely if $H$ is a Boolean algebra.
\begin{proof}
Let $x \# y$ denote the term in the theorem statement. Clearly $x \# y$ satisfies irreflexivity and symmetry in any Heyting algebra. Moreover, $\bot \# y = y$ and $y \# \top = \neg y$ both hold for all $y$. Now assume that $x \# y$ is cotransitive in $H$ as well. Then we can prove that the law of excluded middle holds in $H$ by taking an arbitrary $x \in H$ and arguing as follows:
\begin{align*}
\top &= \bot \# \top                       & & \text{since $\bot \# y = y$ for all $y$}\\
     &\leq (\bot \# x) \vee (x \# \top)    & & \text{by cotransitivity of $\#$}\\
     &= x \vee \neg x                      & & \text{by $\bot \# y = y$ and $y \# \top = \neg y$}
\end{align*}
This shows that if $\#$ is cotransitive then $H$ is a Boolean algebra. Conversely, if $H$ is a Boolean algebra, then an eight-fold case distinction argument in classical propositional logic immediately shows that $x \# y$ satisfies both cotransitivity and tightness.
\end{proof}
\end{proposition}

\begin{proposition}\label{prop:wlem-rijke2}
Rijke's second candidate, the term $\neg (x \leftrightarrow y)$ constitutes an apartness term on a Heyting algebra $H$ precisely if $H$ satisfies the weak law of excluded middle, $\neg x \vee \neg\neg x = \top$ for all $x \in H$.
\begin{proof}
As in Proposition~\ref{prop:lem-rijke1}, the term satisfies irreflexivity and symmetry in any Heyting algebra. And again, it's easy to show that $y \# \top = \neg y$ holds for all $y \in H$. This time, however, we have $\bot \# y = \neg\neg y$ instead of $\bot \# y = y$. Assuming that $x \# y$ is cotransitive in $H$, we can prove that the weak law of excluded middle holds in $H$ by taking an arbitrary $x \in H$ and arguing as follows:
\begin{align*}
\top &= \bot \# \top                       & & \text{since $\bot \# y = y$ for all $y$}\\
     &\leq (\bot \# x) \vee (x \# \top)    & & \text{by cotransitivity of $\#$}\\
     &= \neg\neg x \vee \neg x                      & & \text{by $\bot \# y = \neg\neg y$ and $y \# \top = \neg y$}
\end{align*}
This shows that if $\#$ is cotransitive then $H$ satisfies the weak law of excluded middle.

The converse follows from an easy argument in IPC. Consider two propositions $P,Q$ so that $\neg (P \leftrightarrow Q)$ holds. Observe that $(\neg X \wedge \neg Y) \rightarrow (X \leftrightarrow Y)$ holds by the principle of explosion. Case analysis on weak excluded middle for $P,Q$, combined with this observation, yields that either $\neg P \wedge \neg\neg Q$ or $\neg\neg P \wedge \neg Q$ must hold. From there, a case analysis on weak excluded middle for $R$ quickly yields $\neg (P \leftrightarrow R) \vee \neg (R \leftrightarrow Q)$.
\end{proof}
\end{proposition}

\section{Possible reducts}\label{sec:possible-reducts}

\begin{figure}
    \centering
    \includegraphics[width=0.25\textwidth]{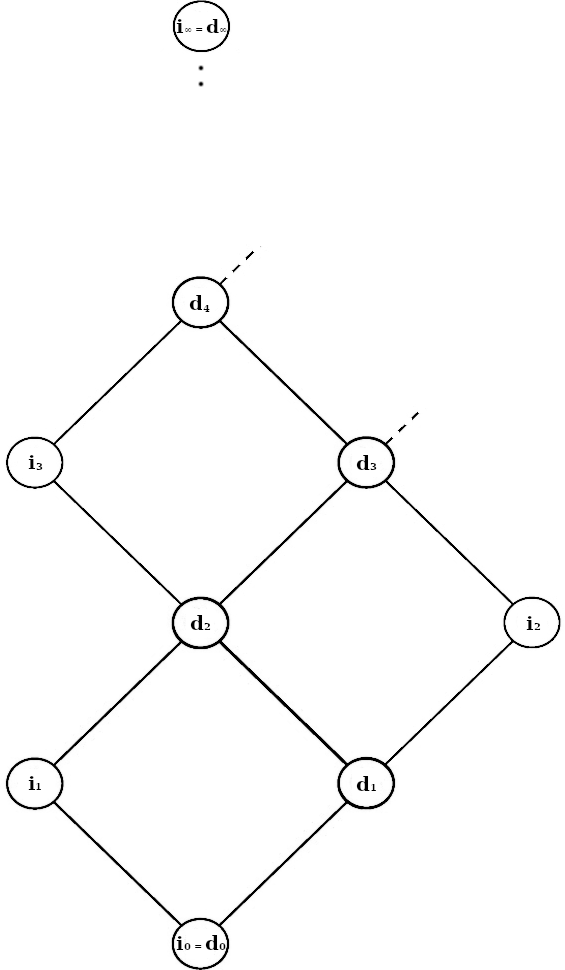}
    \caption{Hasse diagram of the Rieger-Nishimura lattice.}
    \label{fig:rnlattice}
\end{figure}

\begin{definition}\label{def:rieger-nishimura-lattice}
We define the sequences $d$ and $i$ of \emph{disjunctive} and \emph{implicative Rieger-Nishimura terms} by mutual recursion as follows:
\begin{align*}
    & \mathbf{d_0} = \bot, & \quad & \mathbf{i_0} = \bot, \\
    & \mathbf{d_1} = y, & \quad & \mathbf{i_1} = \neg y, \\
    & \mathbf{d_{n+1}} = \mathbf{i_n} \vee \mathbf{d_n}, & \quad & \mathbf{i_{n+1}} = \mathbf{i_n} \rightarrow \mathbf{d_n}, \\
    & \mathbf{d_\infty} = \top & \quad & \mathbf{i_\infty} = \top.
\end{align*}
The \emph{Rieger-Nishimura lattice} consists of the terms $\mathbf{d_n}, \mathbf{i_n}$ for all $n \in \mathbb{N}\cup\{\infty\}$ in the free variable $y$, equipped with the ordering defined by the following clauses:
\begin{enumerate}
    \item $\mathbf{d_n} \leq \mathbf{d_m}$ precisely if $n \leq m$,
    \item $\mathbf{d_n} \leq \mathbf{i_m}$ precisely if $n = 0$ or $n < m$,
    \item $\mathbf{i_n} \leq \mathbf{d_m}$ precisely if $n = 0$ or $n < m$, and
    \item $\mathbf{i_n} \leq \mathbf{i_m}$ precisely if $n = 0$, $n = m$ or $n + 1 < m$.
\end{enumerate}
One can check easily that this construction determines the lattice (and in fact a Heyting algebra) depicted on Figure~\ref{fig:rnlattice}.
\end{definition}

\begin{theorem}[Rieger \cite{rieger-nishimura-lattice}]\label{thm:rieger-nishimura}
The Rieger-Nishimura lattice is isomorphic to the free Heyting algebra on one generator. In particular, the function sending each term of the Rieger-Nishimura lattice to (the logical equivalence class of) the corresponding IPC formula in one variable is an isomorphism of Heyting algebras.
\end{theorem}

\begin{body}
It follows from Theorem~\ref{thm:rieger-nishimura} that every IPC formula in one proposition letter is logically equivalent to a formula corresponding to a term of the Rieger-Nishimura lattice. If $x \# y$ is an apartness term in the Heyting algebra $H$, the term $\top \# y$ corresponds to an IPC formula in one proposition letter, which means that we can find a term $t$ in the Rieger-Nishimura lattice so that $H \models \top \# y = t$. Similarly, we can find a term $b$ in the Rieger-Nishimura lattice so that $H \models \bot \# y = b$. We will classify apartness terms by proving results about the possible values of $t,b$ (the possible \textit{reducts} of $x \# y$ in the terminology of Definition~\ref{def:reducts}).
\end{body}

\begin{definition}\label{def:reducts}
Consider an apartness term $x \# y$ in a Heyting algebra $H$. The \textit{top-reduct} (resp. \textit{bottom-reduct}) of $x \# y$ is the smallest term $t$ in the Rieger-Nishimura lattice so that $H \models \top \# y = t$ (resp. $H \models \bot \# y = t$).
\end{definition}

\begin{lemma}\label{lemma:topreduct}
Take a non-trivial apartness term $x \# y$ in a non-trivial Heyting algebra $H$. The top-reduct of $x \# y$ is $\neg y$.
\begin{proof}
Consider the top-reduct $\top \# y$ of an apartness term $x \# y$ in a Heyting algebra $H$. Since the top-reduct has only one free variable, it coincides with a Rieger-Nishimura term by Theorem~\ref{thm:rieger-nishimura}. We have three cases to consider:
\begin{enumerate}
    \item $\top \# y$ belongs to the upset of $\mathbf{d_1}$ (the term $y$); or
    \item $\top \# y$ is equivalent to $\mathbf{i_0}$ (the term $\bot$); or
    \item $\top \# y$ is equivalent to $\mathbf{i_1}$ (the term $y \rightarrow \bot$).
\end{enumerate}
In the first case, we have
\begin{align*}
\top &= \top \rightarrow (\top \# \top) & & \text{since $y \leq \top \# y$ for all $y$}\\
     &= \top \rightarrow \bot           & & \text{by irreflexivity of $\#$}\\
     &= \bot                            & & \text{by $\neg\top = \bot$}
\end{align*}
so $H$ is the trivial Heyting algebra. In the second case, the apartness term itself is trivial, since for all $x,y \in H$ we have
\begin{align*}
x \# y &= (x \# \top) \vee (\top \# y)  & & \text{by cotransitivity of $\#$}\\
       &= (\top \# x) \vee (\top \# y)  & & \text{by symmetry of $\#$}\\
       &= \bot \vee \bot                & & \text{since $\top \# z = \bot$ for all $z\in H$}\\
       &= \bot                          & & \text{by $\bot \vee \bot = \bot$}
\end{align*}
This leaves only the third case, that the top-reduct of $x \# y$ is $\neg y$.
\end{proof}
\end{lemma}

\begin{theorem}\label{thm:boolean}
A Heyting algebra $H$ that has a tight apartness term is a Boolean algebra.
\begin{proof}
Since $x \# y$ is tight, we get that the equality $\neg (x \# y) \rightarrow (x \rightarrow y) = \top$ holds for all $x,y\in H$. Substituting $\top$ for $x$, we get that $$\neg (\top \# y) \rightarrow (\top \rightarrow y) = \top.$$ From Lemma~\ref{lemma:topreduct} we know that $\top \# y = \neg y$, and since $\top \rightarrow y = y$ in every Heyting algebra, we get that $\neg\neg y \rightarrow y = \top$ holds for all $y \in H$. As noted in the introduction, this means that $H$ is a Boolean algebra with $x \rightarrow y = \neg x \vee y$.
\end{proof}
\end{theorem}

\begin{lemma}\label{lemma:bottomreduct}
Take a non-trivial apartness term $x \# y$ in a non-trivial Heyting algebra $H$. The bottom-reduct $\bot \# y$ is one of either $y$ or else $\neg\neg y$.
\begin{proof}
Consider the bottom-reduct $\bot \# y$ of the term $x \# y$. As usual, by cotransitivity, we cannot have $\bot \# y = \bot$. Since $\bot \# y$ has one free variable, Theorem~\ref{thm:rieger-nishimura} applies, and it suffices to show that $\bot \# y$ cannot belong in the upset of $\neg y$ either. If this happened, we'd have the contradiction
\begin{align*}
\top &= \neg\bot         & & \text{ }\\
     &\leq \bot \# \bot  & & \text{by $\neg y \leq \bot \# y$}\\
     &= \bot             & & \text{by irreflexivity.}
\end{align*}
This shows that the bottom-reduct of $x \# y$ is one of either $y$ or else $\neg\neg y$.
\end{proof}
\end{lemma}

\begin{theorem}\label{thm:characterization}
A Heyting algebra $H$ has a non-trivial apartness term precisely if it satisfies the weak law of excluded middle.
\begin{proof}
We already know (Proposition~\ref{prop:wlem-rijke2}) that Rijke's second candidate is an apartness term in every Heyting algebra which satisfies the weak law of excluded middle.
For the other direction, consider a Heyting algebra $H$ with a non-trivial apartness term $x \# y$. We must show that $H$ satisfies the weak law of excluded middle. By Lemma~\ref{lemma:topreduct}, we know that $\top \# y = \neg y$. By Lemma~\ref{lemma:bottomreduct}, we know that either $\bot \# y = y$ or $\bot \# y = \neg\neg y$. We treat these two cases separately.

\textbf{Case 1}: Assume that $\bot \# y = y$ for all $y \in H$. We will show that in this case, $H$ is a Boolean algebra, and thus satisfies the (weak \textit{and} strong) law of excluded middle. Take an arbitrary $x \in H$ and argue as follows:
\begin{align*}
\top &= \bot \# \top                     & & \text{by $\bot \# y = y$}\\
     &\leq (\bot \# x) \vee (x \# \top)  & & \text{by cotransitivity of $\#$}\\
     &= (\bot \# x) \vee (\top \# x)     & & \text{by symmetry of $\#$}\\
     &= x \vee \neg x                    & & \text{by $\bot \# y = y$ and $\top \# y = \neg y$}
\end{align*}

\textbf{Case 2}: Assume that $\bot \# y = \neg\neg y$ for all $y \in H$. In this case, a nigh-identical argument to that of Case 1 allows us to conclude the weak law of excluded middle. Take an arbitrary $x \in H$, then reason as follows:
\begin{align*}
\top &= \bot \# \top                     & & \text{by $\bot \# y = y$}\\
     &\leq (\bot \# x) \vee (x \# \top)  & & \text{by cotransitivity of $\#$}\\
     &= (\bot \# x) \vee (\top \# x)     & & \text{by symmetry of $\#$}\\
     &= \neg\neg x \vee \neg x           & & \text{by $\bot \# y = \neg\neg y$ and $\top \# y = \neg y$}
\end{align*}
In either case, $H$ satisfies the weak law of excluded middle as claimed.
\end{proof}
\end{theorem}

\begin{corollary}\label{cor:characterization}
A superintuitionistic logic $L$ has a non-trivial apartness relation between propositions precisely if $\neg P \vee \neg\neg P \in L$. In particular, intuitionistic logic has no non-trivial apartness relation between propositions.
\end{corollary}

\section{Classification}\label{sec:classification}

\begin{body}
Lemma~\ref{lemma:bottomreduct} and Theorem~\ref{thm:characterization} allow us to quickly identify necessary conditions for a given formula giving rise to a non-trivial apartness relation in intuitionistic logic. For example, Propositions~\ref{prop:lem-rijke1}~and~\ref{prop:wlem-rijke2} follow immediately by considering the respective bottom-reducts of Rijke's candidates.
\end{body}

\begin{body}
We have already seen that Rijke's first candidate is a non-trivial apartness relation between propositions only in classical propositional logic. Notice that in classical propositional logic, Rijke's second candidate defines the same apartness relation as the first candidate. Here we explain this phenomenon by showing that in a superintuitionistic logic satisfying the weak law of excluded middle, Rijke's second candidate in fact represents the \textit{only} apartness relation between propositions, up to logical equivalence.
\end{body}

\begin{proposition}\label{prop:inequality1}
Take a Heyting algebra $H$ with an apartness term $x \# y$. The inequality $\neg\neg x \wedge \neg\neg y \leq x \# \neg y$ holds for all $x,y \in H$.
\begin{proof}
By Lemma~\ref{lemma:bottomreduct}, we have either $\bot \# y = y$ or $\bot \# y = \neg\neg y$.
Moreover, if $\bot \# y = y$ then $H$ is a Boolean algebra: so in fact we know $\bot \# y = \neg\neg y$ in both cases! Thus, it suffices to show that with the assumptions
\begin{enumerate}
    \item[A1.] $\neg\neg P$,
    \item[A2.] $\neg\neg Q$,
    \item[A3.] $(\bot \# Q) \leftrightarrow \neg\neg Q$,
    \item[A4.] $(\bot \# \neg P) \leftrightarrow \neg\neg\neg P$, and
    \item[A5.] $(\bot \# Q) \rightarrow (\bot \# \neg P) \vee (\neg P \# Q)$
\end{enumerate}
we have $\vdash (\text{A1} \wedge \text{A2} \wedge \text{A3} \wedge \text{A4} \wedge \text{A5}) \rightarrow (P \# \neg Q)$ in intuitionistic logic. The following argument works:
\begin{enumerate}
    \item Assume that A1, A2, A3, A4 and A5 all hold.
    \item Combining A2 and A3 we get that $\bot \# Q$ also holds. 
    \item From $\bot \# Q$, A5 gives $(\bot \# \neg P) \vee (\neg P \# Q)$.
    \item But by A4, $\bot \# \neg P$ is $\neg\neg\neg P$, which contradicts A1. Thus, we must have $\neg P \# Q$. \qedhere
\end{enumerate}
\end{proof}
\end{proposition}

\begin{theorem}\label{thm:classification-lem}
Consider a Boolean algebra $H$. Up to equivalence of terms, Rijke's first candidate is the unique non-trivial apartness term in $H$.
\begin{proof}
We already know from Proposition~\ref{prop:lem-rijke1} that Rijke's first candidate is an apartness term in every Boolean algebra. We must now show that there are no others. So take a Boolean algebra $H$ with a non-trivial apartness term $x \# y$. Using double-negations, Proposition~\ref{prop:inequality1} gives the inequalities $x \wedge \neg y \leq x \# y$ and $\neg x \wedge y \leq x \# y$. Thus, $(x \wedge \neg y) \vee (\neg x \wedge y) \leq x \# y$. For the other direction, notice that the inequalities
\begin{align*}
x \# y &\leq (x \# \bot) \vee (\bot \# y)  & & \text{by cotransitivity of $\#$}\\
       &= (\bot \# x) \vee (\bot \# y)  & & \text{by symmetry of $\#$}\\
       &= x \vee y                      & & \text{by Lemma~\ref{lemma:bottomreduct}}
\end{align*}
and
\begin{align*}
x \# y &\leq (x \# \top) \vee (\top \# y)  & & \text{by cotransitivity of $\#$}\\
       &= (\top \# x) \vee (\top \# y)  & & \text{by symmetry of $\#$}\\
       &= \neg x \vee \neg y                      & & \text{by Lemma~\ref{lemma:topreduct}}
\end{align*}
also hold, so $x \# y \leq (x \wedge \neg y) \vee (\neg x \wedge y)$ follows.
\end{proof}
\end{theorem}

\begin{proposition}\label{prop:inequality2}
In a Heyting algebra $H$ with a non-trivial apartness term $x \# y$, the equality $x \# \neg\neg x = \bot$ holds for all $x \in H$.
\begin{proof}
We have both
\begin{align*}
x \# \neg\neg x &\leq (x \# \top) \vee (\top \# \neg\neg x) & & \text{by cotransitivity of $\#$}\\
                &= (\top \# x) \vee (\top \# \neg\neg x) & & \text{by symmetry of $\#$}\\
                &= \neg x \vee \neg\neg\neg x                & & \text{by Lemma~\ref{lemma:topreduct}}
\end{align*}
and
\begin{align*}
x \# \neg\neg x &\leq (x \# \bot) \vee (\bot \# \neg\neg x) & & \text{by cotransitivity of $\#$}\\
                &= (\bot \# x) \vee (\bot \# \neg\neg x) & & \text{by symmetry of $\#$}\\
                &= \neg \neg x \vee \neg\neg\neg\neg x                & & \text{by Lemma~\ref{lemma:bottomreduct}.}
\end{align*}
This means that $x \# \neg\neg x \leq (\neg x \vee \neg\neg\neg x) \wedge (\neg \neg x \vee \neg\neg\neg\neg x) = \bot$.
\end{proof}
\end{proposition}

\begin{theorem}\label{thm:classification}
Consider a Heyting algebra $H$ satisfying the weak law of excluded middle. Up to equivalence of terms, Rijke's second candidate is the unique non-trivial apartness term in $H$.
\begin{proof}
Take a Heyting algebra $H$ with an arbitrary non-trivial apartness term $x \# y$. We will show the uniqueness of $\#$ by proving the  equality $$x \# y = (\neg\neg x \wedge \neg y) \vee (\neg x \wedge \neg\neg y).$$
Since we already know (Proposition~\ref{prop:wlem-rijke2}) that Rijke's second candidate is an apartness term in $H$, this will immediately imply its equivalence with $x \# y$.

First, notice that the inequalities
\begin{align*}
x \# y &\leq (x \# \bot) \vee (\bot \# y)  & & \text{by cotransitivity of $\#$}\\
       &= (\bot \# x) \vee (\bot \# y)  & & \text{by symmetry of $\#$}\\
       &= \neg\neg x \vee \neg\neg y    & & \text{by Lemma~\ref{lemma:bottomreduct}}
\end{align*}
and
\begin{align*}
x \# y &\leq (x \# \top) \vee (\top \# y)  & & \text{by cotransitivity of $\#$}\\
       &= (\top \# x) \vee (\top \# y)     & & \text{by symmetry of $\#$}\\
       &= \neg x \vee \neg y               & & \text{by Lemma~\ref{lemma:topreduct}}
\end{align*}
both hold, so the inequality $x \# y \leq (\neg\neg x \vee \neg\neg y) \wedge (\neg x \vee \neg y) = (\neg\neg x \wedge \neg y) \vee (\neg x \wedge \neg\neg y)$ holds as well. For the other direction, we can argue as follows:
\begin{align*}
\neg\neg x \wedge \neg y
    &\leq \neg\neg x \wedge \neg\neg\neg y & & \text{by $y \leq \neg\neg y$}\\
    &\leq x \# \neg\neg y                  & & \text{by Proposition~\ref{prop:inequality1}}\\
    &\leq (x \# y) \vee (y \# \neg\neg y)  & & \text{by cotransitivity of $\#$}\\
    &= (x \# y) \vee \bot                  & & \text{by Proposition~\ref{prop:inequality2}}\\
    &= x \# y                              & & \text{by $x \vee \bot = x$}\\
\end{align*}
The same argument gives $\neg x \wedge \neg\neg y = \neg\neg y \wedge \neg x \leq y \# x = x \# y$ as well. Thus, we obtain $(\neg\neg x \wedge \neg y) \vee (\neg x \wedge \neg\neg y) \leq x \# y$ as claimed.
\end{proof}
\end{theorem}

\begin{corollary}
A superintuitionistic logic $L$ has at most one non-trivial apartness relation between propositions. Moreover, any such relation is logically equivalent to Rijke's second candidate.
\end{corollary}

\section{Type Theory}\label{sec:typetheory}

\begin{body}
Possible reduct considerations allow us to neatly classify apartness relations between propositions in superintuitionistic logics. Unfortunately, one cannot push this technique further, to obtain results about apartness relations between propositions in Martin-Löf Type Theory~\cite{martinLofITT}, or even just in second-order propositional logic: Theorem~\ref{thm:rieger-nishimura} has no analogue in these settings.
\end{body}

\begin{body}
Although the possible reducts technique does not generalize, and we certainly cannot have a complete classification of all possible apartness relations that can arise between propositions, parametricity methods do allow us to show at least that, under a propositions-as-types interpretation, plain Martin-Löf Type Theory itself does not define any non-trivial apartness relation between propositions (i.e. between types). In Homotopy Type Theory, where one uses the propositions-as-some-types interpretation, we have a much simpler situation: propositional extensionality allows us to establish the conclusions of Lemma~\ref{lemma:topreduct} and Lemma~\ref{lemma:bottomreduct} directly (as Lemma~\ref{lemma:typed-reducts-converse}), so the proof of Theorem~\ref{thm:characterization} goes through verbatim.
\end{body}

\begin{body}
We work in Martin-Löf Type Theory equipped with a hierarchy of universe types $\mathcal{U} : \mathcal{U}' : \dots$, not assumed to be cumulative. Given a type $P: \mathcal{U}$ and a family $\varphi: P \rightarrow \mathcal{U}$, we write the dependent product (function) type as $\Pi x : P. \varphi(x)$ and the dependent sum (pair) type as $\Sigma x : P. \varphi(x)$. We construct elements of the latter by pairing, writing $(p,q) : \Sigma x : P. \varphi(x)$ given some $p:P$ and $q:\varphi(p)$. We have the standard projections $\mathrm{pr}_1 : (\Sigma x : P. \varphi(x)) \rightarrow P$ and $\mathrm{pr}_2 : \Pi y:(\Sigma x : P. \varphi(x)). \varphi(\mathrm{pr}_1(y))$ obeying the equations $\mathrm{pr}_1(p,q) = p$ and $\mathrm{pr}_2(p,q)=q$. The symbol $\mathbf{0}$ denotes the empty type, $\mathbf{1}$ denotes the set with one element, and $\Bool$ denotes the two-element set with elements $0,1$. Identity types over some type $P$ are written as $\jequiv_P$.
\end{body}

\begin{definition}\label{def:apartness-type}
In Martin-Löf Type Theory, an \textit{apartness relation between types} consists of the following data:
\begin{itemize}
    \item \textbf{Relation}: a term $\# : \mathcal{U} \rightarrow \mathcal{U} \rightarrow \mathcal{U}$,
    \item \textbf{Irreflexive}: a term $\mathrm{irr} : \Pi x: \mathcal{U}. x \# x \rightarrow \mathbf{0}$,
    \item \textbf{Symmetric}: a term $\mathrm{sym} : \Pi x, y: \mathcal{U}. x \# y \rightarrow y \# x$,
    \item \textbf{Cotransitive}: a term $\mathrm{cot} : \Pi x, y, z: \mathcal{U}. x \# y \rightarrow ((x \# z) + (z \# y))$
\end{itemize}
We call the relation \textit{tight} if we also have a term $\mathrm{thg} : \Pi x, y : \mathcal{U}. (x \# y \rightarrow \mathbf{0}) \rightarrow x \jequiv_{\mathcal{U}} y$, and \textit{visibly nontrivial} if we also have a term $\mathrm{vnt} : \Sigma A,B :\mathcal{U}. A \# B$.
\end{definition}

\begin{body}
We now show that Martin-Löf Type Theory does not define any visibly nontrivial apartness relations between types. For this, we use the parametricity results of Bernardy, Jansson and Paterson~\cite{bernardyParametricity}. Recall that parametricity lets us associate a binary relation (a so-called \textit{parametricity conditon}) to each type in such a way that the assumption that every term satisfies the parametricity condition generated by its type can be consistently made in Type Theory.
\end{body}

\begin{proposition}\label{prop:uniformity}
It is consistent with Martin-Löf Type Theory that all functions of signature $\mathcal{U} \rightarrow \Bool$ are constant.
\begin{proof}
The parametricity condition for the type $\mathcal{U} \rightarrow \Bool$ is the following: $$\lambda f_1.\lambda f_2. \Pi A_1,A_2 : \mathcal{U}. \Pi R: (A_1 \rightarrow A_2 \rightarrow \mathcal{U}). (f_1\ A_1) \jequiv_{\Bool} (f_2\ A_2).$$
We can derive this from the definition of parametricity on type schemes of the form $\mathcal{U} \rightarrow B$ for some $B : \mathcal{U}$ (Section~3.3 of~\cite{bernardyParametricity}) and the fact that the identity relation $\jequiv_\Bool$ is the parametricity condition associated with the type $\Bool$ (\textit{ibid}, Section 5.1). The assumption that every term satisfies the parametricity condition generated by its type can be consistently made in Type Theory. But if every function $f: \mathcal{U} \rightarrow \Bool$ satisfies the parametricity condition above, then $\Pi f : \mathcal{U} \rightarrow \Bool. \Pi A_1,A_2 : \mathcal{U}. f\ A_1 \jequiv_{\Bool} f\ A_2$ follows immediately.
\end{proof}
\end{proposition}

\begin{corollary}\label{cor:mlttnondef}
Martin-Löf Type Theory does not define any visibly nontrivial apartness relations between types.
\begin{proof}
Assume it did, with terms $\#, \mathrm{irr}, \mathrm{sym}, \mathrm{cot}, \mathrm{vnt}$. Let $A : \mathcal{U}$ denote $\mathrm{pr}_1 \mathrm{vnt}$, and $B : \mathcal{U}$ denote $\mathrm{pr}_1(\mathrm{pr}_2 \mathrm{vnt})$. Case analysis on the result of $\mathrm{cot}$ lets us define a non-constant function $f:\mathcal{U}\rightarrow \Bool$ which maps $C:\mathcal{U}$ to $0$ if $A\# C$, and to $1$ if $C \# B$. More formally, we can let $f$ denote the function $$\lambda z: \mathcal{U}. \mathrm{rec}_{(A\# z) + (z\# B)}\Bool(\lambda a:(A \# z). 1)(\lambda b: (z \# B). 0)(\mathrm{cot}\ A\ B\ z\ (\mathrm{pr}_2 (\mathrm{pr}_2 \mathrm{vnt})))$$
and inhabit the types $f\ A \jequiv_\Bool 0$ and $f\ B \jequiv_\Bool 1$ as follows. Define
\begin{align*}
  & g_A : \Pi x : (A \# A) + (A \# B). \mathrm{rec}_{(A\# A) + (A\# B)}\Bool (\lambda a.1) (\lambda b.0) x \jequiv_\Bool 0 \\
  & g_A := \mathrm{ind}_{(A \# A) + (A \# B)} (\lambda x. \mathrm{rec}_{(A\# A) + (A\# B)} (\lambda a. 1) (\lambda b. 0) x \jequiv_\Bool 0) (\mathrm{rec}_{\mathbf{0}} \circ \mathrm{irr}\ A) ) (\lambda b. \mathrm{refl}_0) \\
  & g_B : \Pi x : (A \# B) + (B \# B). \mathrm{rec}_{(A\# B) + (B\# B)}\Bool (\lambda a.1) (\lambda b.0) x \jequiv_\Bool 1 \\
  & g_B := \mathrm{ind}_{(A \# B) + (B \# B)} (\lambda x. \mathrm{rec}_{(A\# B) + (B\# B)} (\lambda a. 1) (\lambda b. 0) x \jequiv_\Bool 1) (\lambda a. \mathrm{refl}_1) (\mathrm{rec}_{\mathbf{0}} \circ \mathrm{irr}\ B) )
\end{align*}
to obtain $g_A(\mathrm{cot}\ A\ B\ A\ \kappa) : f\ A \jequiv_\Bool 0$ and $g_B(\mathrm{cot}\ A\ B\ B\  \kappa) : f\ B \jequiv_\Bool 1$ where $\kappa$ abbreviates $\mathrm{pr}_2 (\mathrm{pr}_2 \mathrm{vnt})$. We have managed to construct a function $f: \mathcal{U} \rightarrow \Bool$ that is not constant, contradicting Proposition~\ref{prop:uniformity}.
\end{proof}
\end{corollary}

\begin{body}
Note that for a propositional formula $\varphi$, we have $\vdash \varphi$ precisely if Martin-Löf Type Theory inhabits the propositions-as-types translation of $\varphi$ (see Section 6 of~\cite{awodeyPropositions}). Thus, Corollary~\ref{cor:mlttnondef} yields an alternative proof of the fact that intuitionistic propositional logic has no non-trivial apartness relation between propositions. Similarly, one can obtain an alternative proof of Corollary~\ref{cor:mlttnondef} itself using results about a stronger system. Escardó and Streicher established (Theorem 2.2.~of~\cite{escardoIntrinsic}) that in Homotopy Type Theory, the existence of a function $f: \mathcal{U} \rightarrow \mathbf{2}$ and types $A,B: \mathcal{U}$ satisfying $f\ A \jequiv_\mathcal{U} f\ B \rightarrow \mathbf{0}$ implies an unprovable instance of the law of excluded middle, and hence Homotopy Type Theory cannot construct any function with this property. \textit{A fortiori}, neither can Martin-Löf Type Theory itself. However, as shown above, a visibly nontrivial apartness relation would allow us to construct such an $f$. It again follows that Martin-Löf Type Theory does not define visibly nontrivial apartness relations between types.
\end{body}

\begin{definition}\label{def:escardos-apart}
Let $\mathcal{U}'$ denote the universe one level above $\mathcal{U}$. The relation $\#_2: \mathcal{U} \rightarrow \mathcal{U} \rightarrow \mathcal{U}'$ is defined by the expression $\lambda A. \lambda B. \Sigma p : \mathcal{U} \rightarrow \Bool. p\ A \jequiv_\Bool p\ B \rightarrow \mathbf{0}$.
\end{definition}

\begin{proposition}[M.~Escard\'{o}~\cite{escardoTypeTop}]\label{prop:escardos-apart}
The relation $\#_2$ constitutes an apartness relation between types (modulo universe issues in Definition~\ref{def:apartness-type}, since it targets the universe \textit{above} $\mathcal{U}$).
\end{proposition}

\begin{body}
Escard\'{o}'s results show that Martin-Löf Type Theory does not prove the tightness of the apartness relation $\#_2$. It follows from Proposition~\ref{prop:uniformity} that Martin-Löf Type Theory cannot prove $\#_2$ visibly nontrivial either.
\end{body}

\begin{body}
One possible reading of Theorems~\ref{thm:boolean}~and~\ref{thm:classification} says that the term $$t(x,y) = (\neg\neg x \wedge \neg y) \vee (\neg x \wedge \neg\neg y)$$ acts like a sort of \textit{canary} for the existence of nontrivial and tight apartness between propositions in a given Heyting algebra. I.e., there is a nontrivial apartness between propositions precisely if $t$ is nontrivial, and a tight apartness between propositions precisely if $t$ is tight. The proof of Corollary~\ref{cor:mlttnondef} has a similar reading: it shows that Escard\'{o}'s apartness relation $\#_2$ acts like a canary for visibly nontrivial apartness between types. As always, size/universe issues complicate matters: unlike the IPC canary, this relation takes large values.
\end{body}

\begin{proposition}
Martin-Löf Type Theory has a canary for visibly nontrivial tight apartness.
\begin{proof}
Let $\mathcal{U}'$ denote the universe one level above $\mathcal{U}$. Let $\mathrm{IsTap}\ R$ abbreviate the data associated with a tight apartness in Definition~\ref{def:apartness-type} for the relation $R: \mathcal{U} \rightarrow \mathcal{U} \rightarrow \mathcal{U}$. Define the relation $\#_\star: \mathcal{U} \rightarrow \mathcal{U} \rightarrow \mathcal{U}'$ as $$\lambda A. \lambda B. \Sigma R : \mathcal{U} \rightarrow \mathcal{U} \rightarrow \mathcal{U}. (\mathrm{IsTap}\ R) \times (R\ A\ B).$$
As an exercise, the reader can show that $\#_\star$ is always an apartness relation. Now assume that there is a visibly nontrivial tight apartness $\#$. Let $A : \mathcal{U}$ denote $\mathrm{pr}_1 \mathrm{vnt}$, and $B : \mathcal{U}$ denote $\mathrm{pr}_2 \mathrm{vnt}$. We have $\left(\#, \mathrm{pr}_2 (\mathrm{pr}_2 \mathrm{vnt})\right) : A \#_\star B$,
showing that $\#_\star$ is also visibily nontrivial. For tightness, we sketch the construction of the relevant term. Have $P,Q : \mathcal{U}$ and assume that $a: P \#_\star Q \rightarrow \mathrm{\mathbf{0}}$. Assume for a contradiction that $b: P\# Q$ holds. Then we would have $(\#, b): P \#_\star Q$, and thus $a(\#, b) : \mathrm{\mathbf{0}}$. Discharging the assumption $b$, we have $\lambda b. a(\#,b) : P\# Q \rightarrow \mathbf{0}$. Therefore, $\mathrm{tgh}\ P\ Q\ (\lambda b. a(\#,b)) : P \jequiv_\mathcal{U} Q$. This shows that if there is a visibly nontrivial tight apartness, then in particular $\#_\star$ is a visibly nontrivial tight apartness.
\end{proof}
\end{proposition}

When a type admits a tight apartness relation, its equality is stable under double-negation. In particular:

\begin{theorem}\label{thm:stability}
In Martin-Löf Type Theory with a tight apartness, the type $P \jequiv_\mathcal{U} Q$ is always stable under double-negation.
\begin{proof}[Proof sketch]
Prove $(P \# Q) \rightarrow ((P \jequiv_\mathcal{U} Q) \rightarrow \mathbf{0})$. Taking contrapositives, we get $(((P \jequiv_\mathcal{U} Q) \rightarrow \mathbf{0}) \rightarrow \mathbf{0}) \rightarrow ((P \# Q) \rightarrow \mathbf{0})$. But by tightness we already know that the implication $((P \# Q) \rightarrow \mathbf{0}) \rightarrow (P \jequiv_\mathcal{U} Q)$ holds, so we're done.
\end{proof}
\end{theorem}

\subsection{Propositions-as-some-types}

\begin{body}
In the propositions-as-some-types paradigm, we identify propositions with subsingleton types: a subsingleton type $X$ is one so that $\Pi x:X. \Pi y:X. x \jequiv_{X} y$ is inhabited. Some authors call this \textit{univalent logic}, in contrast with the \textit{Curry--Howard logic} of propositions-as-types, emphasizing the relevance of this paradigm to univalent foundations (see e.g. slide 9 of~\cite{escardoSlides}). For instance, the propositions-as-types version of the law of excluded middle contradicts the univalence axiom (Corollary~3.2.7 of \cite{hottbook}), whereas the propositions-as-some-types version does not (slides 27 and 28 of~\cite{escardoSlides}).
\end{body}

\begin{body}
In the remainder of this section, we study the behavior of apartness relations between propositions within the \textit{propositions-as-some-types} paradigm. We show that one can prove analogues of Lemma~\ref{lemma:topreduct} and Lemma~\ref{lemma:bottomreduct} directly (without using the reducts technique) under a propositions-as-some-types interpretation, allowing us to carry out the arguments of Sections~\ref{sec:possible-reducts}~and~\ref{sec:classification} directly in Type Theory.
\end{body}

\paragraph{Conventions}

\begin{body}
We write $\Omega$ for the type of (subsingleton) propositions $\Sigma P:\mathcal{U}. \Pi x:P. \Pi y:P. x \jequiv_{U} y$. By a slight abuse of notation, when this causes no ambiguity, we use the same name $X$ to refer to both $X : \Omega$ and the corresponding $\mathrm{pr}_1\ X: \mathcal{U}$.
\end{body}

\begin{body}
Assumptions such as the univalence axiom and bracket/truncation types usually accompany propositions-as-some-types. To carry out the proofs in this section, it suffices to assume that we have access to terms $\pet, \peb$ of types
\begin{align*}
  & \pet : \Pi P:\Omega. P \rightarrow (\mathbf{1} \jequiv_\mathcal{U} P) \\
  & \peb : \Pi P:\Omega. (P \rightarrow \mathbf{0}) \rightarrow (P \jequiv_\mathcal{U} \mathbf{0}).
\end{align*}
Keep in mind that these both follow from the principle of \textit{propositional extensionality},
$$ \mathrm{pe} : \Pi P,Q : \Omega. (P \rightarrow Q) \rightarrow (Q \rightarrow P) \rightarrow (P \jequiv_{\mathcal{U}} Q).$$
In turn, propositional extensionality follows as an immediate consequence of the univalence axiom of Homotopy Type Theory. Therefore, all our results apply to Homotopy Type Theory and related systems such as Cubical Type Theory~\cite{cohenCuTT} verbatim.
\end{body}

\begin{body}
Throughout the section, we assume that apartness relations between subsingletons have the signature $\# : \Omega \rightarrow \Omega \rightarrow \mathcal{U}$, and that cotransitivity takes the form $$\mathrm{cot} : \Pi x, y, z: \Omega. x \# y \rightarrow ((x \# z) + (z \# y)),$$ just as it did when we dealt with propositions-as-types. However, in type theories that come equipped with propositional truncations, all constructions of this section carry over to strictly subsingleton-valued relations satisfying ``mere cotransitivity'', i.e. to signatures $\# : \Omega \rightarrow \Omega \rightarrow \Omega$ and $$\mathrm{cot} : \Pi x, y, z: \Omega. x \# y \rightarrow \|(x \# z) + (z \# y)\|.$$
\end{body}

\begin{lemma}\label{lemma:typed-reducts}
In a type theory with $\pet \mathrm{/} \peb$, every irreflexive binary relation $\#$ between subsingletons satisfies $(\mathbf{1} \# P) \rightarrow (P \rightarrow \mathbf{0})$ and $(\mathbf{0} \# P) \rightarrow ((P \rightarrow \mathbf{0}) \rightarrow \mathbf{0})$.
\begin{proof}
Along with $\# : \Omega \rightarrow \Omega \rightarrow \mathcal{U}$ and $\mathrm{irr} : \Pi x: \Omega. x \# x \rightarrow \mathbf{0}$, take $P : \Omega$, $q: \mathbf{1} \# P$ and $p : P$. We want to construct a term of type $\mathbf{0}$ using these. We have $\pet\ P\ p : \mathbf{1} \jequiv_\mathcal{U} P$, which allows us to write
\begin{align*}
  & g: \mathbf{1} \# \mathbf{1} \\
  & g:= \mathrm{transport}^{\lambda x. \mathbf{1}\# x} (\pet\ P\ p)\ q
\end{align*}
and therefore $\mathrm{irr}\ P\ g : \mathbf{0}$ as desired. A similar construction works for the $\mathbf{0} \# P$ case.
\end{proof}
\end{lemma}

\begin{body}
Our goal is to derive analogues of Lemma~\ref{lemma:topreduct} and Lemma~\ref{lemma:bottomreduct} in type theory. Notice that Lemma~\ref{lemma:typed-reducts} already gives us implications between the apartness terms and their purported reducts in one direction: the converse direction has a few prerequisites, and will follow in Lemma~\ref{lemma:typed-reducts-converse}.
\end{body}

\begin{theorem}\label{thm:types-boolean}
In a type theory with propositional extensionality, the existence of a tight, irreflexive relation between subsingletons implies double-negation elimination.
\begin{proof}
We shall imitate the proof of Theorem~\ref{thm:boolean}, using Lemma~\ref{lemma:typed-reducts} in place of Lemma~\ref{lemma:topreduct}. Start with
\begin{align*}
  & \# : \Omega \rightarrow \Omega \rightarrow \mathcal{U} \\
  & \mathrm{irr} : \Pi x: \Omega. x \# x \rightarrow \mathbf{0} \\
  & \mathrm{thg} : \Pi x, y : \Omega. (x \# y \rightarrow \mathbf{0}) \rightarrow x \jequiv_{\mathcal{U}} y
\end{align*}
along with $P: \Omega$ and $p: (P \rightarrow \mathbf{0}) \rightarrow \mathbf{0}$. We seek to write down something of type $P$. Let $g : \mathbf{1} \# P \rightarrow (P \rightarrow \mathbf{0})$ denote the term obtained from Lemma~\ref{lemma:typed-reducts}. We then have that
$p \circ g : \mathbf{1} \# P \rightarrow \mathbf{0}$, and so we can set
\begin{align*}
  & h : \mathbf{1} \jequiv_\mathcal{U} P \\
  & h := \mathrm{thg}\ \mathbf{1}\ P\ (p \circ g)
\end{align*}
and obtain $\mathrm{transport}^{\lambda x. x}\ h\ \star : P$ as desired.
\end{proof}
\end{theorem}

\begin{body}
Alternatively, one can establish Theorem~\ref{thm:types-boolean} using Theorem~\ref{thm:stability}, by first inhabiting $\neg\neg P \rightarrow \neg\neg (\mathbf{1} \jequiv_\mathcal{U} P)$ using the double contrapositive of $\mathrm{pe}_\mathbf{1}$, then invoking $\neg\neg (\mathbf{1} \jequiv_\mathcal{U} P) \rightarrow \mathbf{1} \jequiv_\mathcal{U} P$ and $\mathrm{transport}$ to conclude $\neg\neg P \rightarrow P$.
\end{body}

\begin{proposition}\label{prop:nontriviality}
In a type theory with propositional extensionality, a visibly nontrivial irreflexive, symmetric, cotransitive relation $\#$ between subsingletons always satisfies $\mathbf{1} \# \mathbf{0}$.
\begin{proof}
Consider such a relation $\#$ with terms $\mathrm{irr}, \mathrm{sym}, \mathrm{cot}, \mathrm{vnt}$. Let $A : \mathcal{U}$ denote $\mathrm{pr}_1 \mathrm{vnt}$, and $B : \mathcal{U}$ denote $\mathrm{pr}_1(\mathrm{pr}_2 \mathrm{vnt})$. For each $P: \Omega$ let $g_P : \mathbf{1} \# P \rightarrow P \rightarrow \mathbf{0}$ denote the term obtained from Lemma~\ref{lemma:typed-reducts}. Construct
\begin{align*}
  & h_P : (\mathbf{1} \# P) \rightarrow (\mathbf{1} \# \mathbf{0}) \\
  & h_P := \lambda t. \mathrm{transport}^{\lambda x. \mathbf{1} \# x} (\peb\ P\ (g_P\ t))\ t \\
  & j : (A \# \mathbf{1}) + (\mathbf{1} \# B) \\
  & j := \mathrm{cot}\ A\ B\ \mathbf{1}\ (\mathrm{pr}_2(\mathrm{pr}_2 \mathrm{vnt}))
\end{align*}
and conclude with $\mathrm{rec}_{(A \# \mathbf{1}) + (\mathbf{1} \# B)} (\mathbf{1} \# \mathbf{0}) (h_A \circ \mathrm{sym}\ A\ \mathbf{1})\ h_B\ j : \mathbf{1} \# \mathbf{0}$ as required.
\end{proof}
\end{proposition}

\begin{lemma}\label{lemma:typed-reducts-converse}
In a type theory with $\pet / \peb$, a visibly nontrivial irreflexive, symmetric, cotransitive relation $\#$ between subsingletons always satisfies $(P \rightarrow \mathbf{0}) \rightarrow (\mathbf{1} \# P)$ and $((P \rightarrow \mathbf{0}) \rightarrow \mathbf{0}) \rightarrow (\mathbf{0} \# P)$.
\begin{proof}
We prove the latter - the former follows by a very similar construction. Consider such a relation $\#$ with terms $\mathrm{irr}, \mathrm{sym}, \mathrm{cot}$, and (by Proposition~\ref{prop:nontriviality}) a term $\mathrm{vnt} : \mathbf{1} \# \mathbf{0}$. Given $P : \Omega$ and $p : (P \rightarrow \mathbf{0}) \rightarrow \mathbf{0}$ we need to construct something of type $\mathbf{0} \# P$. Let $g : \mathbf{1} \# P \rightarrow P \rightarrow \mathbf{0}$ denote the term obtained from Lemma~\ref{lemma:typed-reducts}. Then we have
\begin{align*}
  & j : (\mathbf{1} \# P) + (P \# \mathbf{0}) \\
  & j := \mathrm{cot}\ \mathbf{1}\ \mathbf{0}\ P\ \mathrm{vnt}
\end{align*}
and so $\mathrm{rec}_{(\mathbf{1} \# P) + (P \# \mathbf{0})}\ (\mathbf{0} \# P)\ (\mathrm{rec}_\mathbf{0}\ (\mathbf{0} \# P) \circ p \circ g)\ (\mathrm{sym}\ P\ \mathbf{0})\ j : \mathbf{0} \# P$ as required.
\end{proof}
\end{lemma}

\begin{body}
Using Lemma~\ref{lemma:typed-reducts} and Lemma~\ref{lemma:typed-reducts-converse}, we get that all results of Sections~\ref{sec:possible-reducts}~and~\ref{sec:classification} go through unchanged for type theory, in the propositions-as-some-types setting. In particular, up to logical equivalence, Rijke's second candidate constitutes the unique candidate for non-trivial apartness relations between subsingletons.
\end{body}

\subsection*{Future Work}

\begin{body}
Many interesting questions remain regarding apartness relations between types. Can we find other canaries for the existence of visibly nontrivial apartness relations between types, say ones that do not ``step outside'' the universe $\mathcal{U}$? And what is the logical strength of the existence of a visibly nontrivial / tight apartness relation between types? In particular, is there any term in the Rieger-Nishimura lattice whose universal closure (under the propositions-as-types paradigm) is implied by the existence of a visibly nontrivial apartness? We suspect that the answer is negative: if we do not need to preserve canonicity, it seems like it should be possible to extend type theory consistently with a ``generic'' (even tight) apartness relation without affecting the propositional fragment.
\end{body}

\begin{body}
Finally, note that in our study of apartness relations between subsingletons, we effectively translated meta-theoretic findings from intuitionistic logic into theorems of propositions-as-some-types Type Theory by leveraging the principles $\pet/\peb$. This recurring, emerging pattern, which often allows for deriving results \textit{in} Type Theory from external proof-theoretic arguments (see~\cite{zoltanProofTheory} for an application), remains intriguingly elusive: further investigation and explicitation of it is certainly warranted.
\end{body}

\subsection*{Acknowledgements}

We thank the anonymous referees for many valuable comments and suggestions which helped improve the quality of the manuscript, including the informal summary of the proof strategy in Corollary~\ref{cor:mlttnondef} and the alternative proof of the same corollary, as well as the alternative proof of Theorem~\ref{thm:types-boolean} using Theorem~\ref{thm:stability}.

\bibliography{biblio}
\bibliographystyle{plainurl}

\end{document}